\documentclass{amsart}

\input{cyracc.def}
\newfont{\cyrfnt}{wncyr10}
\newcommand{\cyr}{\baselineskip12.pt\cyrfnt\cyracc}

\usepackage{amssymb,enumerate,verbatim}

\numberwithin{equation}{section}

\theoremstyle{plain}
\newtheorem{theorem}{Theorem}[section]
\newtheorem{corollary}[theorem]{Corollary}
\newtheorem{lemma}[theorem]{Lemma}

\theoremstyle{definition}

\newtheorem{remark}[theorem]{Remark}
\newtheorem{note}[theorem]{Note}
\newtheorem{definition}[theorem]{Definition}
\newtheorem{example}[theorem]{Example}

\newcommand{\Lan}{L}

\newcommand{\Lans}{L^*}

\newcommand{\C}{\mathcal{C}}

\newcommand{\Cl}{\mathcal{K}}

\newcommand{\Z}{\mathbb{Z}}

\newcommand{\N}{\mathbb{N}}

\newcommand{\K}{\mathbb{C}}

\newcommand{\acl}{\operatorname{acl}}

\newcommand{\acls}{\operatorname{acl}^*}

\newcommand{\dcl}{\operatorname{dcl}}

\newcommand{\cl}{\operatorname{cl}}

\newcommand{\trdeg}{\operatorname{trdeg}}

\newcommand{\Mr}{\operatorname{MR}}
\newcommand{\Mrs}{\operatorname{MR}^*}

\newcommand{\tp}{\operatorname{tp}}
\newcommand{\tps}{\operatorname{tp}^*}

\newcommand{\sym}{\operatorname{Sym}}
\newcommand{\dd}{\operatorname{d}}

\begin{document}

\title{On fields and colors}

\date{January 14, 2005}
\author{A. Baudisch}
\address{Institut f\"ur Mathematik\\Humboldt-Universit\"at zu Berlin\\ D-10099 Berlin, Germany}
\email{baudisch@mathematik.hu-berlin.de}

\author{A. Martin-Pizarro}
\address{Institut f\"ur Mathematik\\Humboldt-Universit\"at zu Berlin\\ D-10099 Berlin, Germany}
\email{pizarro@mathematik.hu-berlin.de}

\author{M. Ziegler}
\address{Mathematisches Institut\\
         Albert-Ludwigs-Universit\"at Freiburg\\ D-79104 Freiburg, Germany}
\email{ziegler@uni-freiburg.de}

\thanks{\noindent A. Martin-Pizarro conducted research supported by
the DFG under Forschungsstipendium MA3310/1-1 \\ \indent {\cyr
Professor Martin Cigler yavlyaet\cydot sya postoyannym gostem
universiteta imeni Gumbol\cprime dta v Berline.}}  


\begin{abstract}
We exhibit a simplified version of the construction of a field of
Morley rank $p$ with a predicate of rank $p-1$, extracting the main
ideas for the construction from previous papers and refining the
arguments. Moreover, an explicit axiomatization is given, and ranks are computed.
\end{abstract}

\maketitle

\section{Introduction}\label{S:intro}

\mbox{}

Zil'ber posed the question whether or not every strongly minimal set
whose geo\-metry was not locally modular arose from an algebraic curve
over an algebraically closed field. The conjecture, true in the case
of Zariski Geometries \cite{HZ93}, remained open until E. Hrushovski
\cite{Hr93} refuted it developing a procedure, taking ideas from
Fra\"\i ss\'e, in order to construct countable structures with a
richer and more complicated geometry starting from simpler
ones. Moreover, he was able to merge two algebraically closed fields
of different characteristics into one strongly minimal set
\cite{Hr92}. This procedure was later adapted by Poizat \cite{Po99} to
obtain an algebraically closed field of any given characteristic with
a predicate (whose elements were called \emph{black}, after some
considerations on the \emph{political correctness} of such a choice of
terminology) such that the field has Morley rank $\omega 2$ and the
black points $\omega$. He then used Hrushovski's collapsing method and
produced ``rich'' fields of rank $2$ with black points of rank $1$,
provided the rich field is $\omega$--saturated. A proof of
$\omega$--saturation was supplied by Baldwin and Holland
(\cite{BH00}). Poizat and Baldwin \& Holland also explained how to
obtain fields of rank $p$ with a predicate of rank $1$ and $p-1$,
respectively.

The main goal of this work is to give a complete self-contained proof
of the above facts simplifying as much as possible the arguments. One
of the novelties of this work is exhibiting an axiomatization for the
resulting theory, obtained by direct translation of Hrushovski's
fusion article \cite{Hr92} to the case of colored fields. Actually, we
use a simplified version (see \cite{BMZ04}) of the aforementioned
article, following the spirit of Poizat's black points.

All throughout this work a saturated
enough algebraically closed field $\K$ of some given characteristic $q$ and a
natural number $p\geq 2$. We will prove the following:
\newtheorem*{teo}{Main Theorem}
\begin{teo}[\cite{Po99},\cite{BH00}]
  $\K$ has a subset $N$ such that $(\K,N)$ has Morley rank $p$ and $N$
  has Morley rank $p-1$.
\end{teo}

This paper is structured as follows: We first consider finite partial substructures of $\K$ with some points colored in black. A $\delta$ function is introduced, and Hrushovski's codes \cite{Hr92} are used to described minimal extension (with a small correction from their original definition). The number of certain such extensions is bounded with a $\mu$ function. In this case, we can proceed with the collapse, and the resulting structure is a rich field as in \cite{Po99}. We show that rich fields are exactly the $\omega$-saturated models of a given theory, whose axioms are explictly given. Finally, we compute the Morley rank in terms of $\delta$.
 
This work originated from a seminar held at Humboldt-Universit\"at zu
Berlin directed by the first and third author \cite{BZ04} during
2003--2004 in which the second author took part. We would specially
like to thank Juan R. Bueno for his help in a preliminary version of
this article \cite{Bu04} during his stay in Berlin.\\

\section{Codes}
\newcounter{code_enum}
In this section, we work exclusively inside $\K$. All formulae
are $\Lan$--formulas, where $\Lan$ is the ring language
\begin{definition}\label{D:code}
A \emph{code} $\alpha$ is a tuple
consisting of the following objects: Natural numbers $n_\alpha$,
$m_\alpha$, $k_\alpha$ and formulae
$\varphi_\alpha(\vec{x},\vec{y})$ and
$\psi_\alpha(\vec{x}_1,\dotsc,\vec{x}_{m_\alpha},\vec{y})$ such that
the following holds (We will write $\theta_\alpha(\vec{y})=\exists
\vec{x}\,\varphi_\alpha(\vec{x},\vec{y})$):

\begin{enumerate}[(i)]
\item $\mathrm{length}(\vec{x})=\mathrm{length}(\vec{x_i})=n_\alpha$
\item\label{D:code_grad} If $\models \theta_\alpha(\vec{b})$, then
  $\varphi_\alpha(\vec{x},\vec{b})$ has Morley rank $k_\alpha$ and
  degree $1$.
\item\label{D:code_algebraisch} Let $\vec{a}\models
\varphi_\alpha(\vec{x},\vec{b})$ be generic. For $s\subset
\{1,\dotsc,n_\alpha\}$, write $a_s=\{a_j\}_{j\in s}$. Then, for every
$i\leq n_\alpha$ and $\vec{a}'\models
\varphi_\alpha(\vec{x},\vec{b}')$, we have that:
$$\begin{array}{rcl}
a_i \in \acl(a_s,\vec{b})  &\Longrightarrow&  a'_i \in \acl(a'_s,\vec{b}')\\
a_i \in a_s\vec{b}  &\Longrightarrow&  a'_i \in a'_s\vec{b}'\\
a_i \not\in a_s\vec{b}  &\Longrightarrow&  a'_i \not\in a'_s\vec{b}'
\end{array}$$
\item\label{D:code_eindeutig} If $\models\theta_\alpha(\vec{b})$, then
  $\;\;\Mr\bigl(\varphi_\alpha(\vec{x},\vec{b})\vartriangle\,
  \varphi_\alpha(\vec{x},\vec{b}')\bigr)<k_\alpha\;\Longrightarrow\;
  \vec{b}=\vec{b}'$.
\item\label{D:code_dcl}
  $\models \psi_\alpha(\vec{x}_1,\dotsc,\vec{x}_{m_\alpha},\vec{b})$
  implies that
  $\vec{b}\in\dcl(\vec{x}_1,\dotsc,\vec{x}_{m_\alpha})$.\footnote{%
  Note that the definable closure $\dcl(\vec{x})$ is the perfect hull
  of the field generated by $\vec{x}$.}
\item $\psi_\alpha(\vec{x}_1,\dotsc,\vec{x}_{m_\alpha},\vec{b})$ is
  consistent for all $\vec{b}\models\theta_\alpha(\vec{y})$.
\item\label{D:code_erweiterung} Given
$(\vec{a}_1,\dotsc,\vec{a}_{m_\alpha},\vec{b})\models \psi_\alpha$ and
a generic $\vec{a}'$ realizing $\varphi_\alpha(\vec{x}, \vec{b})$ and
independent from $\vec{a}_1,\ldots,\vec{a}_{m_\alpha}$ over $\vec{b}$, it
follows that
\[\bigwedge\limits_{i=1}^{m_\alpha}
\psi_\alpha(\vec{a}_1,\vec{a}_2,\dotsc,\vec{a}_{i-1},\vec{a}',
\vec{a}_{i+1},\dotsc,\vec{a}_{m_\alpha},\vec{b})\]
\noindent holds.
\setcounter{code_enum}{\value{enumi}}
\end{enumerate}
\end{definition}

\begin{lemma}\label{L:cb}
If $\models\theta_\alpha(\vec{b})$, we have that $\vec{b}$ is a
canonical basis of the type of Morley rank $k_\alpha$ determined by
$\varphi_\alpha(\vec{x},\vec{b})$. 
\end{lemma}
\begin{proof} 
This follows immediately from (\ref{D:code_eindeutig}). 
\end{proof}

\begin{lemma} 
For each definable set $X$ of $n$-tuples of Morley rank $k$ and degree
$1$ there is a code $\alpha$ with $n_\alpha=n$, $k_\alpha=k$
and some $\vec{b}\models \theta_\alpha(\vec{y})$ such that $ \Mr
\bigl( \varphi_\alpha(\vec{x},\vec{b})\vartriangle X \bigr)<k_\alpha$.
\end{lemma} 
\begin{proof}
  Let $X$ be given. We begin with a formula
  $\varphi(\vec{x},\vec{b}_0)$ such that $\vec{b}_0$ is a canonical base
  of the type determined by $X$ and such that $\Mr \bigl(
  \varphi(\vec{x},\vec{b}_0)\vartriangle X \bigr)<k$. Since Morley rank
  and degree are definable in algebraically closed fields, we may assume
  that $\varphi$ satisfies (\ref{D:code_grad}). If, in addition,
  $\varphi(\vec{x},\vec{b}_0)$ witnesses all algebraic dependencies and
  equalities between the components of a generic solution, property
  (\ref{D:code_algebraisch}) holds also. Now, $\vec{b_0}$ is a canonical
  base if and only if $\Mr\bigl(\varphi(\vec{x},\vec{b}_0) \vartriangle
  \varphi(\vec{x},\vec{b}') \bigr)<k\;\to\;\vec{b}_0=\vec{b}'$ for every
  $\vec{b}'\models\tp(\vec{b}_0)$. Thus $\varphi$ satisfies
  (\ref{D:code_eindeutig}), if we add some finite part of
  $q(\vec{y})=\tp(\vec{b}_0)$ to $\varphi(\vec{x},\vec{y})$.
  
  Choose generic
  realizations $\vec{a}_1,\dotsc,\vec{a}_m$ of
  $\varphi(\vec{x},\vec{b}_0)$, independent over $\vec{b}_0$. If $m$ is
  large enough, we have $\vec{b}_0=f(\vec{a}_1,\dotsc,\vec{a}_m)$ for
  some $0$--definable function $f$. We strengthen $\varphi$, so that
  $\vec{b}=f(\vec{x}_1,\dotsc,\vec{x}_m)$ for every sequence
  $\vec{x}_1,\dotsc,\vec{x}_m$ of independent generic realizations of
$\varphi(\vec{x},\vec{b})$, and set $\varphi_\alpha=\varphi$.
  
Finally let
$\psi_\alpha(\vec{x}_1,\dotsc,\vec{x}_m,\vec{b})$ express\\

\parbox{32em}{``For every sequence $\vec{x}_{m+1},\dotsc,\vec{x}_{2m}$
of generic realizations of $\varphi_\alpha(\vec{x},\vec{b})$, such
that $\vec{x}_{1},\ldots,\vec{x}_{2m}$ is independent over $\vec{b}$,
and every choice of distinct indices
$i_1,\dotsc,i_m\in\{1,\dotsc,2m\}$, we have
$\vec{b}=f(\vec{x}_{i_1},\dotsc,\vec{x}_{i_m})$''.\qedhere}
\end{proof}

Let $\alpha$ be a code and $\sigma$ a permutation of
$\{1,\dotsc,n_\alpha\}$. We denote by $\alpha^\sigma$ the code
obtained from $\alpha$ by permuting each of the tuples $\vec{x}$,
$\vec{x}_1$,\ldots,$\vec{x}_{m_\alpha}$ in $\varphi_\alpha$ and
$\psi_\alpha$ according to $\sigma$. We call $\alpha^\sigma$ a
\emph{permutation} of $\alpha$.

\begin{definition}
  Two codes $\alpha$ and $\alpha'$ are equivalent if
  $n_\alpha=n_{\alpha'}$ and $m_\alpha=m_{\alpha'}$ and
  \begin{itemize}
  \item for every realization $\vec{b}$ of $\theta_\alpha$ there is a
    tuple $\vec{b}'$ such that (in $\K$)
    $\varphi_\alpha(\vec{x},\vec{b})\equiv
    \varphi_{\alpha'}(\vec{x},\vec{b}')$ and
    $\psi_\alpha(\vec{x}_1,\dots,\vec{x}_{m_\alpha},\vec{b})\equiv
    \psi_{\alpha'}(\vec{x}_1,\dots,\vec{x}_{m_\alpha},\vec{b}')$.
    \item the same replacing the roles of $\alpha$ and $\alpha'$.
  \end{itemize}
\end{definition}

The following lemma is a slightly weaker as the statement of Lemma 2 in
\cite{Hr92}.
\begin{lemma}
  There is a set $\C$ of codes such that
  \begin{enumerate}[\upshape(i)]
    \setcounter{enumi}{\value{code_enum}}
  \item\label{L:code_cexistenz} For each (non--empty) definable set
    $X$ of Morley degree $1$ there is a code $\alpha\in\C$ and some
    $\vec{b}$ such that\/ $ \Mr \bigl(
    \varphi_\alpha(\vec{x},\vec{b})\vartriangle X \bigr)<k_\alpha$.
  \item\label{L:code_ceindeutig} If $\alpha,\alpha'\in \C$,
    $\models\theta_\alpha(\vec{b})$ and
    $\Mr\bigl(\varphi_\alpha(\vec{x},\vec{b})
    \vartriangle\varphi_{\alpha'}(\vec{x},\vec{b}') \bigr)<k_\alpha$,
    then $\alpha'=\alpha$.\footnote{We identify two codes if there
    defining formulas are equivalent in $\K$.}
  \item\label{L:code_permutation} If $\alpha$ belongs to $\C$, then
    each permutation of $\alpha$ is equivalent to a code in $\C$.
  \end{enumerate}
\end{lemma}
\begin{proof}
  We refer to the claim of (\ref{L:code_cexistenz}) as ``$X$ can be
  coded by $\alpha$''.  List all
  non--empty definable sets of degree $1$ up to conjugation by
  automorphisms of $\K$ by $X_1,X_2,\dotsc$ This is possible since
  $\textrm{ACF}_q$ is small, i.e. it has only countably many $n$-types
  for each $n$. It is enough to show that each $X_i$ can be coded by
  some elements of $\C$.  We will obtain $\C$ as the union of a sequence
  $\emptyset=C_0\subset C_1\subset\dotsb$ of finite sets of codes,
  constructed as follows. Assume that $C_{i-1}$ has been
  constructed and it is closed under permutations in the
  weak sense of (\ref{L:code_permutation}). If $X_i$ can be coded by
  an element of $C_{i-1}$, we set $C_i=C_{i-1}$. Otherwise, choose a
  code $\alpha$ and $\vec{b}_0$ such that $\Mr \bigl(
  \varphi_\alpha(\vec{x},\vec{b}_0)\vartriangle X \bigr)<k_\alpha$. We
  replace $\varphi_\alpha$ by
  $$\varphi_\alpha(\vec{x},\vec{y})\land\textrm{``$\{\vec{x}\,|\,
  \varphi_\alpha(\vec{x},\vec{y})\}$ cannot be coded by an element of
  $C_{i-1}$''}.$$ and obtain a new code, which still codes
  $X_i$. We may assume that
  no permutation of $\alpha$ can code a set which can also
  be coded by a code in $C_{i-1}$. Let $G$ be the group of all
  $\sigma\in\sym(n_\alpha)$ such that 
  $$\Mr \bigl( \varphi_\alpha(\vec{x},\vec{b}_0)\vartriangle
  \varphi_{\alpha^\sigma}(\vec{x},\sigma\vec{b}_0)\bigr)<k_\alpha$$
  for some element denoted as $\sigma\vec{b}_0$ which has the same type as $\vec{b}_0$.
  After adding a finite part of the type of $\vec{b}_0$ to
  $\varphi_\alpha(\vec{x},\vec{y})$ we may assume that for all
  realizations $\vec{b}$ of $\theta_\alpha$ and all $\sigma$, there
  exists $\sigma\vec{b}$ with $\Mr \bigl(
  \varphi_\alpha(\vec{x},\vec{b})\vartriangle
  \varphi_{\alpha^\sigma}(\vec{x},\sigma\vec{b})\bigr)<k_\alpha$ iff
  $\sigma\in G$. Note that $\sigma\vec{b}$ is a $\emptyset$-definable
  function of $\vec{b}$. If we let permutations act on the right on
  codes, this defines a left action of $G$ on $\theta_\alpha(\K)$.

  It is easy to check that $$\varphi_\beta(\vec{x},\vec{y})=
  \bigwedge_{\sigma\in G}\varphi_{\alpha^\sigma}(\vec{x},\sigma\vec{y})$$
  and $\psi_\beta(\vec{x}_1,\dotsc,\vec{y})= \bigwedge_{\sigma\in
  G}\psi_{\alpha^\sigma}(\vec{x}_1,\dotsc,\sigma\vec{y})$ defines a
  code, which again codes $X$. Also, for $\sigma\in G$, we have
  $\varphi_\beta(\vec{x},\vec{y})\equiv
  \varphi_{\beta^\sigma}(\vec{x},\sigma\vec{y})$ and
  $\psi_\beta(\vec{x}_1,\dotsc,\vec{y})\equiv
  \psi_{\beta^\sigma}(\vec{x}_1,\dotsc,\sigma\vec{y})$, which shows
  that $\beta$ is equivalent to $\beta^\sigma$. Now choose
  representatives $\rho_1,\dots,\rho_r$ for the right cosets of $G$ in
  $\sym(n_\alpha)$ and set
  $C_i=C_{i-1}\cup\{\beta^{\rho_1},\dotsc,\beta^{\rho_r}\}$.
\end{proof}

\begin{remark} 
Note that the proof holds in a more general setting of a countable
strongly minimal theory with the DMP (\emph{definable multiplicity
property}) where imaginary parameters $\vec{b}$ are allowed. It is not
possible to find $\C$ closed under permutations (as stated in \cite{Hr92}).
\end{remark}

\section{$\delta$-nonsense}\label{S:delta}

Let $X$ be a set. A function $\delta: \mathcal{P}_\mathrm{fin}(X) \to
\Z$ is a \emph{$\delta$-function} if it satisfies the following:
\begin{enumerate}
\item\label{C:empty} $\delta(\emptyset)=0$
\item\label{C:lmod} $\delta(A\cup B)+\delta(A\cap B)\leq \delta(A)+\delta(B)$
\end{enumerate}
Moreover, if for all $A$ we have that $\delta(A)\geq 0$, then we say
that $\delta$ is \emph{nonnegative}.

For finite subsets $A$ and $B$, we define the \emph{relative} $\delta$-value of $A$ over $B$ by:
\[\delta(A/B)=\delta(A\cup B)-\delta(B)\]
Now, (\ref{C:lmod}) is equivalent to $\delta(A/B)\leq \delta(A/A\cap B)$. It is easy to see that for any $A\cap B\subset C\subset B$, we have that $\delta(A/B)\leq \delta(A/C)$.

Hence, we can extend the definition of the relative $\delta$ to subsets $Y$ (possibly not finite) as follows:
\[ \delta(A/Y)=\underset{A\cap Y\subset C\subset Y}{\inf} \delta(A/C)\]

Note that $\delta(A/Y)$ is in $\{-\infty\}\cup\Z$. Using notation from
\cite{Hr93}, we say that $Y$ is \emph{self-sufficient} in $X$ (denoted
as $Y\leq X$) if for all finite $A\subset X$, we have that
$\delta(A/Y)\geq 0$. We have that
$\delta$ is nonnegative iff $\emptyset\leq X$.

$Y$ is self-sufficient iff $\delta(A)\geq \delta(A\cap Y)$ for all
$A$. If $Y\leq X$, it follows that $Y\cap Z\leq
Z$ for all $Z$. Hence, self-sufficiency is transitive. Moreover, the
intersection of self-sufficient sets is again self-sufficient and each
set $S$ is contained in a smallest self-sufficient subset, its
\emph{self-sufficient closure} $\cl_X(S)$. If $\delta$ is nonnegative,
finite sets have finite closures.

A proper extension $Y\leq Z$ is \emph{minimal} if no $Y\subsetneq
Y'\subsetneq Z$ is self-sufficient in $Z$.  The extension $Z\setminus Y$ must be
finite, which allows us to express minimality by \[\delta(Z/Y')<0 \text{ for all } Y\subsetneq Y'\subsetneq Z\].

\section{Black points}
We extend the ring language $\Lan$ to $\Lans=\Lan\cup \{N\}$, where
$N$ is a unary predicate. \textbf{All considered $\Lans$--structures
are \emph{colored} subsets of $\K$}, i.e. subsets $A$ of $\K$ endowed
with an interpretation $N(A)$ for $N$ (\emph{les points noirs}).
The notation $A\subset B$ implies $N(A)=A\cap N(B)$.

We want to amalgamate \`a la Fra\"\i ss\'e-Hrushovski finite
$\Lans$--structures $A$ according to a function $\delta$ defined 
as follows:

\begin{center}
$\delta(A)=p\cdot\trdeg A-|N(A)|$
\end{center}

Note that $\delta$ satisfies conditions (\ref{C:empty}) and (\ref{C:lmod})
from Section \ref{S:delta}. With this particular definition, we have
that $\delta(\{a\})\leq p$. We are in a setting as in the previous
section. 

Although the general amalgam was studied in careful detail in
\cite{Po99}, we will concentrate on the collapse closer to the spirit
of \cite{Hr92}. Hence, we will consider just sets, and not the
$\Lans$-substructures that they
generate. Nonetheless, in an abuse of notation, we
will call them $\Lans$-structures (and not \emph{partial}
$\Lans$-structures).\\

All the lemmas in the rest of the section are true for arbitrary,
finite or infinite, $\Lans$--structures.
\begin{lemma}\label{L:minext}
Let $B\leq A$ be a minimal extension. We have one of the following cases:
\begin{enumerate}
\item\label{white} If $A$ contains a white point $a$ not in $B$, then $A=B\cup\{a\}$. Moreover, $\delta(A/B)=0$ or $p$, depending whether $a$ is algebraic or transcendental over $B$.
\item\label{black} Otherwise, $A=B\cup\{a_1,\dotsc,a_n\}$ with
$a_1,\dotsc,a_n$ distinct black and $0\leq \delta(A/B)\leq
p-1$. Moreover, for any $\emptyset\neq S\subsetneq
\{a_1,\dotsc,a_n\}$, we have that
$$p\cdot\trdeg(A/B\cup S)< n-|S|.$$
 If $\delta(A/B)= p-1$, then $A=B\cup\{a\}$ with $a$ transcendental over $B$ and black. 
\end{enumerate}
\end{lemma}
\begin{proof} Recall that $B\leq A$ is minimal if it is proper and for any $B\subsetneq A'\subsetneq A$, we have that $\delta(A/BA')<0$. Equivalently, $\delta(A/B)$ is the minimum among all values of $\delta(A'/B)$, where $B\subsetneq A'\subset A$, and it is attained only at $A$.

If $a$ in $A\setminus B$ is white, case (\ref{white}) follows, since we have that $\delta(A\setminus\{a\}/B)\leq \delta(A/B)$, hence $A=B\cup\{a\}$.  The two possibilities for $\delta(A/B)$ are now clear.

Let us assume that $A\setminus B$ contains no white point. Take some $a\in A\setminus B$. Since $B\leq A$, it follows that $a$ is transcendental over $B$ and $\delta(a/B)=p-1$. By minimality, $\delta(A/B)\leq p-1$. 

If $\delta(A/B)=p-1$, then clearly $A=B\cup\{a\}$. 

\end{proof}

\begin{definition}\label{D:good} A minimal extension $B\leq A$ of type 
(\ref{black}) is \emph{good} if $\tp(A/B)$ is
stationary and $\delta(A/B)=0$.
A code $\alpha$ is \emph{good} if it is ``the code'' of a
good minimal extension. That is,

\begin{itemize} 
\item $n_\alpha=p k_\alpha$.
\item $\varphi_\alpha(\vec{x},\vec{y})$ implies that all $x_i$'s are different and different from the components of $\vec{y}$.
\item If $\models\varphi_\alpha(\vec{a},\vec{b})$, for each $\emptyset\neq s\subsetneq\{1,\dotsc,n_\alpha\}$, we have
$$p\cdot\trdeg(\vec{a}/\vec{a}_s\vec{b})<(n_\alpha-|s|).$$ 
\end{itemize}
\end{definition}

\noindent Note that, by (\ref{D:code_algebraisch}), the last two
conditions are true, if they hold for just one realization $\vec{b}$
of $\theta_\alpha$ and one generic realization $\vec{a}$ of
$\varphi_\alpha(\vec{x},\vec{b})$.\\

\noindent Let $\C_g$ be the subset of good codes in $\C$.\\

The next lemma is clear from the definitions.
\begin{lemma}\label{L:code-gut}
  Let $\alpha$ be a good code, $\vec{b}\in\dcl(B)$ realize $\theta_\alpha$, and $\vec{a}$ be a $B$--generic black realization of $\varphi_\alpha(\vec{x},\vec{b})$. Then, $B\cup\{a_1,\dots,a_{n_\alpha}\}$ is a good extension of $B$.
\end{lemma}

\begin{lemma}\label{L:gute-code}
  Let $B\leq A=B\cup\{a_1,\dotsc,a_n\}$ be a good extension. Then
  there is a good code $\alpha$ and $\vec{b}\in\dcl(B)$ such that
  $\models\varphi_\alpha(\vec{a},\vec{b})$.
\end{lemma}
\begin{proof}
  Choose $\chi(\vec{x})\in\tp(\vec{a}/B)$ of Morley rank
  $k=\Mr(\vec{a}/B)$ and degree $1$. There is $\alpha\in\C$ and
  $\vec{b}$ such that
  $\Mr\bigl(\chi(\vec{x})\vartriangle\varphi_\alpha(\vec{x},
  \vec{b})\bigr)<k=k_\alpha$.  Since $\vec{a}$ is a $B$--generic
  realization of $\chi(\vec{x})$, it is also a $B$--generic
  realization of $\varphi_\alpha(\vec{x},\vec{b})$. Since $\vec{b}$ is
  a canonical base of $\tp(\vec{a}/B)$, $\vec{b}$ belongs to
  $\dcl(B)$. Since $A/B$ is good, we have that $\alpha$ is a good code.
\end{proof}
In the previous Lemma, we chose $\vec{a}$ as a $B$--generic realization of
$\varphi_\alpha(\vec{x},\vec{b})$. The following result shows that this the only possibility.
\begin{lemma}[cf.\ Lemma 3A in \cite{Hr92}]\label{L:aux}
Let $\alpha$ be a good code, $\vec{b}\in\acl(B)$ realize $\theta_\alpha$, and $\vec{a}$ be a black realization of $\varphi_\alpha(\vec{x},\vec{b})$ which does not completely lie in $B$.  Then, the following holds:
\begin{enumerate}
\item\label{C:1} $\delta(\vec{a}/B)\leq 0$
\item\label{C:2} If $\delta(\vec{a}/B)= 0$, then $\vec{a}\cap
B=\emptyset$ and $\vec{a}$ is a $B$--generic realization of $\varphi_\alpha(\vec{x},\vec{b})$.
\end{enumerate}
\end{lemma}
\begin{proof}
  If $\vec{a}$ is not
  disjoint from $B$ we have
  $$\delta(\vec{a}/B)\leq\delta(\vec{a}/\vec{a}'\vec{b})<0.$$ for
  $\vec{a}\cap B=\vec{a}'$. Hence, $\delta(\vec{a}/B)=0$ yields that  $\vec{a}\cap B=\emptyset$. In this case, we have $\delta(\vec{a}/B)\leq p\cdot
  k_\alpha-n_\alpha=0$. Therefore, $\trdeg(\vec{a}/B)=k_\alpha$ and $\vec{a}$ is a $B$--generic
  solution of $\varphi_\alpha(\vec{x},\vec{b})$.
\end{proof}

\section{The (in)famous $\mu$ function}\label{S:class}
  We now fix a function
$\mu^*:\C_g\to\N$ which is finite-to-one on the
set of all $\alpha$ with $n_\alpha=n$ for each $n$ in $\N$. Moreover, $\mu^*(\beta)=\mu^*(\alpha)$ must hold if $\beta$ is equivalent to a
permutation of $\alpha$, and
$$\mu^*(\alpha)\geq m_\alpha-1.$$
The function $\mu$ is then defined by 
$$\mu(\alpha)=((p-1)(n_\alpha-1) +1) m_\alpha+\mu^*(\alpha).$$ 
We note that $\mu(\alpha)\geq m_\alpha$.
\begin{note}\label{fte-to-one}
One can replace in the following $\mu$ by
$\mu'(\alpha)=F(\alpha)+\mu^*(\alpha)$ for any function $F$ which
satisfies $F(\alpha^\sigma)=F(\alpha)$ and
$F(\alpha)\geq((p-1)(n_\alpha-1) +1) m_\alpha$. The class of functions
$\mu$ is not increased by this, only the complete theories $T^\mu$ (see Section \ref{S:axioms}) get
weaker, but equivalent, axiomatizations.
\end{note}
We recover the definition introduced in \cite{Po99} for approximations to a Morley sequence of a given good minimal extension. 
\begin{definition}\label{D:PMS}
Let $\alpha$ be a good code and $\vec{b}\models \theta_\alpha$. A
\emph{pseudo-morley sequence} for $\alpha$ over $\vec{b}$ is a
(finite) sequence $\vec{a}_1,\dotsc,\vec{a}_r$ of disjoint
realizations of $\varphi_\alpha(\vec{x},\vec{b})$ painted in black
such that any distinct $m_\alpha$ elements among
$\{\vec{a}_1,\dotsc,\vec{a}_r\}$ realize
$\psi_\alpha(\vec{x}_1,\dotsc,\vec{x}_{m_\alpha},\vec{b})$.
\end{definition}
It follows that $\vec{b}$ is in the definable closure of the pseudo-morley
sequence if $r\geq
m_\alpha$ from part (\ref{D:code_dcl}) of \ref{D:code}.

We now consider the class of $\Lans$-structures on which $\delta$ is
non-negative and for any good code in $\C$, we cannot find a pseudo-morley
sequence that is longer than the value of $\mu$ at this code.

\begin{definition}\label{D:class}
The class $\Cl^\mu$ is the class of all $\Lans$-structures $M$ (i.e\ colored
subsets of $\K$) such that:
\begin{itemize}
\item $\emptyset\leq M$.
\item No $\alpha$ in $\C_g$ has a pseudo-morley sequence in $M$ of
  length longer than $\mu(\alpha)$.
\end{itemize}
\end{definition}
\noindent We denote by $\Cl_\mathrm{fin}^\mu$ the class of all finite
 $\Lans$-structures in $\Cl^\mu$.

Recall that the first condition means that for any finite set
$A\subset M$, we have $\delta(A)\geq 0$. Clearly,
$\Cl_\mathrm{fin}^\mu$ is not empty ($\emptyset$ is an element of this
class). In fact all finite subsets of $\K$ with no black points are
in the class.

Since $\textrm{ACF}_q$ is small, $\Cl_\mathrm{fin}^\mu$ contains at most countably many structures up to isomorphism.
 
The following result resumes the ingredients used in \cite{Hr92}
stating them in a form closer to the original idea of Fra\"\i ss\'e's
amalgamation procedure to construct a coun\-table ultrahomogeneous
model whose age is exactly $\Cl_\mathrm{fin}^\mu$. Moreover, it yields
explicit conditions for an $\Lans$-structure to be a member of
$\Cl^\mu$, which will be useful for exhibiting an axiomatization of
this class.

\begin{lemma}\label{L:member}
Let $M$ be in $\Cl^\mu$ and $M\leq M'$ a minimal extension. 

If $M'$ contains a new white point, then $M'$ is in $\Cl^\mu$. 

Otherwise, $M'$ is in $\Cl^\mu$ if and only if none of the following two
conditions holds:
\begin{enumerate}
\item[a)] There is a code $\alpha\in\C_g$ and a realization $\vec{b}\in
  \dcl(M)$ of $\theta_\alpha$, such that:
  \begin{enumerate}
  \item[i)] $M'\setminus M$ contains a realization
    $\vec{a}$ of $\varphi_\alpha(\vec{x},\vec{b})$.
  \item[ii)] $M$ contains a pseudo-morley sequence for $\alpha$ over
    $\vec{b}$ of length $\mu(\alpha)$.
  \end{enumerate}

\item[b)] There is some code $\alpha\in\C_g$ and a pseudo-morley
sequence for $\alpha$ in $M'$ of length $\mu(\alpha)+1$, such that
there are more than $\mu^*(\alpha)$ many elements of the sequence
contained in $M'\setminus M$.
\end{enumerate}
If a) holds, $\vec{a}$ is an enumeration of $M'\setminus M$ and an
$M$--generic realization of $\varphi_\alpha(\vec{x},\vec{b})$.
\end{lemma}
Since $n_\alpha\geq p$ for good codes, the lemma implies that
$M'\in\Cl^\mu$ if $\delta(M'/M)=p-1$.

\begin{proof} If $M'$ contains a new white point $a$, 
by \ref{L:minext} (\ref{white}), we get that $M'=M\cup\{a\}$. If $M'$
is not in $\Cl^\mu$, it contains a pseudo-morley sequence of length
$\mu(\alpha)+1$ for some code $\alpha\in\C_g$. Since $\{a\}$ adds no
black points, the sequence is contained in $M$, which is a
contradiction.

Suppose now that $M'\setminus M$ has no new white points.  If b)
holds, $M'$ is not in $\Cl^\mu$ by definition. If we have case a),
$\vec{a}$ is generic over $M$ by Lemma \ref{L:aux} (\ref{C:2})
and we can extend the sequence of ii) by $\vec{a}$, thanks to
condition (\ref{D:code_erweiterung}). This shows that
$M'$ is not in
$\Cl^\mu$. Also, since $M'/M$ is minimal and $\delta(\vec{a}/M)=0$, we
have that $M'=M\cup\{a_1,\dotsc,a_{n_\alpha}\}$.

For the other direction, if $M'$ as above is not in $\Cl^\mu$, there
exists a code $\alpha\in\C_g$ and a pseudo-morley sequence
$\vec{e}_0,\dotsc,\vec{e}_{\mu(\alpha)}$ for $\alpha$ in $M'$ over
some $\vec{b}\in \dcl(M')$. We may rearrange the sequence as follows:
\begin{itemize} 
\item $\vec{e}_0,\dotsc,\vec{e}_{r_0-1}$ are contained in $M$.
\item $\vec{e}_{r_0},\dotsc,\vec{e}_{r_1-1}$ are not in
  $M$, but have at least one
  coordinate in $M$.
\item $\vec{e}_{r_1},\dotsc,\vec{e}_{\mu(\alpha)}$ are in $\{a_1,\dotsc,a_{n_\alpha}\}$.
\end{itemize}
Since $M$ is in $\Cl^\mu$, we have that $r_0\leq\mu(\alpha)$. There
are two possibilities:\\

\noindent Case 1. $m_\alpha\leq r_0$. In this case,
$\vec{b}\in\dcl(\vec{e}_0,\dotsc,\vec{e}_{m_\alpha-1})$ is in $\dcl(M)$. By
\ref{L:aux} (\ref{C:1}), we have that $\delta(\vec{e}_{r_0}/M)\leq 0$. Since
$M\leq M'$, we have that $\delta(\vec{e}_{r_0}/M)= 0$. Hence, for each
$i$, we conclude from \ref{L:aux} (\ref{C:2}) that either $\vec{e}_i$ is
disjoint from or contained in $M$. That is, $r_0=r_1$. As above, we conclude
$M'=M\cup\vec{e}_{r_0}$. Hence, $r_0=\mu(\alpha)$ by disjointness of the pseudo-morley sequence. Therefore, a) holds.\\

\noindent Case 2. $r_0\leq m_\alpha$. Define
  $\delta(i)=\delta(\vec{e}_i/M\vec{e}_0,\dotsc,\vec{e}_{i-1})$. Then, since
$M\leq M'$, we have that:
\[0 \leq  \delta(\vec{e}_0,\dotsc,\vec{e}_{r_1-1}/M)=\sum\limits_{i<r_1} \delta(i)=\sum\limits_{i<m_\alpha} \delta(i) +\sum\limits_{m_\alpha\leq i<r_1} \delta(i)\]

For $i<m_\alpha$, we have that $\delta(i)\leq (p-1)(n_\alpha-1)$ (Note that if $\vec{d}$ is a tuple of black points, we always have that $\delta(\vec{d}/B)\leq (p-1)\cdot\trdeg(\vec{d}/B)$ for any set $B$).

For $m_\alpha\leq i<r_1$, it follows that $\vec{b}$ belongs to
$\dcl(M\vec{e}_0\dots\vec{e}_{i-1})$. But there is some
coordinate of $\vec{e}_i$ in $M$, and hence, again from Lemma
\ref{L:aux} (\ref{C:2}), we conclude that $\delta(i)<0$.

From the inequalities above, we get:
\[0\leq (p-1)(n_\alpha-1)m_\alpha-(r_1-m_\alpha)\]
That is, $r_1\leq ((p-1)(n_\alpha-1)+1)m_\alpha$. Now, 
\[\mbox{}\hspace{1cm}\mu(\alpha)-r_1+1\geq \mu(\alpha)-((p-1)(n_\alpha-1)+1) m_\alpha+1\geq\mu^*(\alpha)+1\] 
This yields b).
\end{proof}

\begin{corollary}\label{C:member}
Let $M$ be in $\Cl^\mu$, $\alpha\in\C_g$, $\vec{b}\in\dcl{M}$ a
realization of $\theta_\alpha$ and $\vec{a}$ a black $M$--generic
realization of $\varphi_\alpha(\vec{x},\vec{b})$. Then
$M'=M\cup\{a_1,\dotsc,a_{n_\alpha}\}$ is in $\Cl^\mu$ if and only if
none of the following two conditions holds:
\begin{enumerate}
\item[a)] $M$ contains a pseudo-morley sequence for $\alpha$ over
    $\vec{b}$ of length $\mu(\alpha)$.
\item[b)] There is some code $\beta\in\C_g$ and a pseudo-morley sequence
for $\beta$ in $M'$ of length $\mu(\beta)+1$, such that there are more
than $\mu^*(\beta)$ many elements of the sequence contained in
$M'\setminus M$.
\end{enumerate}
\end{corollary}
\begin{proof}
  $M'$ is a minimal extension of $M$ by \ref{L:code-gut}, so we can
  apply the last lemma. We need only show the following: If
  $\alpha'$ is a code in $\C_g$, $\vec{b}'\in\dcl(M)$ such that $\vec{a}$
  is a permuted $M$-generic realization of
  $\varphi_{\alpha'}(\vec{x},\vec{b}')$ and if $\alpha'$ has a
  pseudo-morley sequence of length $\mu(\alpha')$ in $M$ over
  $\vec{b}'$, then $\alpha$ has a pseudo-morley sequence of length
  $\mu(\alpha)$ in $M$ over $\vec{b}$.

  Let $\sigma$ be a permutation of $\alpha$ such that $\vec{a}$
  realizes $\varphi_{\alpha^\sigma}$. By (\ref{L:code_permutation}) there
  is a code $\alpha''\in\C_g$ which is equivalent to $\alpha^\sigma$. So
  there is $\vec{b}''$ such that
  $\varphi_{\alpha^\sigma}(\vec{x},\vec{b}')\equiv
  \varphi_{\alpha''}(\vec{x},\vec{b}'')$ and
  $\psi_{\alpha^\sigma}(\vec{x}_1,\dotsc,\vec{b}')\equiv
  \psi_{\alpha''}(\vec{x}_1,\dotsc,\vec{b}'')$. The permuted
  pseudo-morley sequence of $\alpha'$ is a pseudo-morley sequence of
  $\alpha''$ over $\vec{b}''\in\dcl(M)$, and $\vec{a}$ is an
  $M$-generic realization of $\varphi_{\alpha''}(\vec{x},\vec{b}'')$.
  The properties (\ref{L:code_ceindeutig}) and
  (\ref{D:code_eindeutig}) of $\C$ imply $\alpha''=\alpha$ and
  $\vec{b}''=\vec{b}$. Finally, we have
  $\mu(\alpha)=\mu(\alpha'')=\mu(\alpha^\sigma)$.
\end{proof}

\section{Fra\"\i ss\'e limits for $\Cl^\mu$}\label{S:Amal}

In this section, we show that the class $\Cl^\mu$ (and hence,
$\Cl_\mathrm{fin}^\mu$) has the Amalgamation Property, and hence, we
can obtain \emph{rich} fields as introduced by Poizat in \cite{Po99}
(We apologize for translating notation into other languages).

An \emph{isomorphism} between two colored subsets $A$ and $B$ of $\K$
is a bijection which maps $N(A)$ onto $N(B)$ and is elementary as a
partial map defined on $\K$. A \emph{self-sufficient embedding} from
$A$ to $B$ is an isomorphism between $A$ and a self-sufficient subset
of $B$.
 
\begin{theorem}\label{T:Amal}
The class $\Cl^\mu$ has the amalgamation property with respect to
self-sufficient embeddings.
\end{theorem}
\begin{proof}
  Let $B\leq M$ and $B\leq A$ be structures in $\Cl^\mu$. We need
  to show that there is an extension $M'$ of $M$ in
  $\Cl^\mu$, with $M\leq M'$ and some $B\leq A'\leq M'$ such that $A$
  and $A'$ are isomorphic over $B$.  By splitting the extension $B\leq
  A$ into minimal ones, we may assume it is minimal.\\
  
  \noindent Case 1. $B\leq A$ has a new white point $a$. Let $p$ be the
  type of $a$ over $B$.  We distinguish two (non-exclusive) cases.\\
    
  \noindent Subcase 1.1. $p$ is algebraic and realized in $M$, say by
  $a'$. Self-sufficiency of $B$ in $M$ yields that $a'$ is white. So
  $A'=B\cup\{a'\}$ is isomorphic to $A$. Since $\delta(a'/B)=0$, it implies
  that $A'\leq M$.\\
  
  \noindent Subcase 1.2.  $p$ can be realized in an extension of $M$ by
  a new element $a'$. We paint $a'$ white and set $M'=B\cup\{a'\}$.\\
  
  \noindent Case 2. $B\leq A$ has no new white points. Since $B$ is
  self-sufficient in $A$, no element of $A\setminus B$ is algebraic
  over $B$.  So we can take for $M'$ be the \emph{free amalgam} (as in
  \cite{Po99}) of $M$ and $A$ over $B$, that is, we assume $M$ and $A$
  to be algebraically independent over $B$ and let $M'$ be their
  union. It is easy to see that $M$ and $A$ are self-sufficient in
  $M'$ and that $M'/M$ is minimal. We are done if $M'$ belongs to
  $\Cl^\mu$. Otherwise, by Lemma \ref{L:member}, there are two cases:\\
  
  \noindent Subcase 2.1) There is a code $\alpha\in\C_g$, a realization
  $\vec{b}$ of $\theta_\alpha$ in $\dcl(M)$, a pseudo-morley
  sequence for $\alpha$ in $M$ over $\vec{b}$ of length
  $\mu(\alpha)$ and $M'\setminus M=\vec{a}$ is a $M$--generic
  realization of $\varphi_\alpha(\vec{x},\vec{b})$. Since $\vec{a}$
  is independent from $M$ over $B$ and $\vec{b}$ is the canonical
  parameter of $\tp(\vec{a}/M)$), we have that $\vec{b}\in\acl(B)$ . The
  sequence cannot be contained in $B$, since $A$ is in
  $\Cl^\mu$. Hence, there is some
  \emph{black} realization $\vec{a}'$ of
  $\varphi_\alpha(\vec{x},\vec{b})$ in $M$ not completely contained
  in $B$. Now, $\delta(\vec{a}'/B)=0$ since $B\leq M$, therefore
  $\vec{a}'$ is generic over $B$ by Lemma \ref{L:aux}. So,
  $B\cup\{a'_1,\dotsc,a'_{n_\alpha}\}$ is self-sufficient in $M$
  and it is isomorphic to $A$
  over $B$.\\
  
  \noindent Subcase 2.2) There is a code $\alpha\in\C_g$, a canonical
  basis $\vec{b}$ in $\dcl(M')$ for $\alpha$ such that there is a
  pseudo-morley sequence $\vec{e}_0,\dotsc,\vec{e}_{\mu(\alpha)}$ for
  $\alpha$ over $\vec{b}$ in $M'$ with more than $\mu^*(\alpha)$ many
  elements coming from $M'\setminus M$. Again, since
  $\mu^*(\alpha)+1\geq m_\alpha$, we have that $\vec{b}$ is in
  $\dcl(A)$. There must be at least one member $\vec{e}_i$ not
  contained in $A$ (because $A$ is in $\Cl^\mu$). Since $A\leq M'$, it
  follows from \ref{L:aux} (\ref{C:2}) that $\vec{e}_i$ is an
  $A$--generic realization of $\varphi_\alpha(\vec{x},\vec{b})$ in
  $M$. But $\vec{e}_i$ and $A$ are independent over $B$, therefore the
  canonical basis $\vec{b}$ of $\alpha$ is in $\acl(B)$. \\
  \noindent
  Pick some $\vec{e}_j$ in $M'\setminus M=A\setminus
  B$. Again it follows that $\vec{e}_j$ is a $B$--generic realization
  of $\varphi_\alpha(\vec{x},\vec{b})$.  Since $M'/M$ is minimal,
  $\vec{e}_j$ enumerates $M'\setminus M$ and we are in subcase
  2.1. Note that all $e_k$, $k\not=j$, are in $M$, which implies
  $\vec{b}\in\dcl(M)$.
\end{proof}

We call $M$ in $\Cl^\mu$ \emph{rich} if for any $B\leq M$ finite and
any finite extension $B\leq A$ of members of $\Cl^\mu$, there is a
self-sufficient substructure $A'\leq M$ with $B\leq A'$ and
$B$-isomorphic to $A$.

\begin{corollary}\label{C:fraisse}
  There is a unique $($up to isomorphism$)$ countable rich structure
  $M$ in $\Cl^\mu$.
\end{corollary}
We will see in Theorem \ref{T:sat=rich} that rich structures are
colored algebraically closed fields. We will call them \emph{rich fields}.

\begin{remark}\label{R:dimcode}
Let $M$ be a rich field\footnote{The remark is true for all models
of $T^\mu$, as defined in Section \ref{S:axioms}.}, $\alpha$ be a code in $\C_g$ and $\vec{b}$ be a realization of
$\theta_\alpha$ in $M$. Let $\dim_\alpha(M/\vec{b})$ be the maximal
length of a pseudo-morley sequence of $\alpha$ over $\vec{b}$ in
$M$ and $B=\cl_M(\vec{b})$ the (finite) self-sufficient closure of
$\vec{b}$ in $M$.
Then there are two cases, either
$$\dim_\alpha(M/\vec{b})=\dim_\alpha(B/\vec{b})$$
or
$$\dim_\alpha(M/\vec{b})=\mu(\alpha).$$
\end{remark}
\begin{proof}
  Choose a black $B$--generic realization $\vec{a}$ of
  $\varphi_\alpha(\vec{x},\vec{b})$ outside $M$. If $A=B\cup\vec{a}$
  is not in $\Cl^\mu$, all realizations of
  $\varphi_\alpha(\vec{x},\vec{b})$ in $M$ are contained in $B$ (by
  Lemma \ref{L:aux} (\ref{C:2})). Therefore
  $\dim_\alpha(M/\vec{b})=\dim_\alpha(B/\vec{b})$.

  If $A$ belongs to $\Cl^\mu$, let $C\leq M$ be a finite extension of
  $B$ with $\dim_\alpha(M/\vec{b})=\dim_\alpha(C/\vec{b})$ and let
  $C'$ be the free amalgam of $C$ and $A$ over $B$. Since $M$
  is rich, $C'$ does not belong to $\Cl^\mu$. The proof of
  \ref{T:Amal} (applied to $C$ instead of $M$) and of \ref{C:member}
  shows that $\dim_\alpha(C/\vec{b})=\mu(\alpha)$.
\end{proof}

\section{A theory for $\Cl^\mu$}\label{S:axioms}

In this section, we will show that the class $\Cl^\mu$ is
axiomatizable and we will give explicit axioms that describe some
completion. Rich fields will then be $\omega$-saturated models of this
theory. First, a foreword about the choice of axioms:

We will see in Section \ref{S:ranks} that extensions with $\delta=0$
will become algebraic. We know (by reducing it to the case of good
minimal extensions) that at most there are $\mu$ many realizations. If
we are given a minimal extension $B\leq A$, where $B\leq M$, we could
amalgamate $A$ and $M$ freely over $B$ and the amalgam could be
potentially an element of $\Cl^\mu$. By richness, this cannot happen,
since there is one realization too many in the amalgam not in
$M$. Hence, we need to prohibit the amalgam to be an element of
$\Cl^\mu$. We know exactly by \ref{L:member} and \ref{C:member}
when this happens. Therefore, our axioms should state that such
an amalgam cannot happen.

The theory $T^\mu$ in the extended language $\Lans=\Lan\cup\{N\}$
 has the following axioms (more precisely, \emph{axiom schemes}):

\begin{description}
\item[Universal Axioms]\mbox{}
\begin{enumerate}
\item\label{intdom} Any model is an integral domain of characteristic $q$.
\item\label{ss} $\emptyset$ is self-sufficient in any model of $T^\mu$.
\item\label{nosol} Given a code $\alpha\in\C_g$, any pseudo-morley
sequence for $\alpha$ has length at most $\mu(\alpha)$.
\end{enumerate}
\item[$\forall\exists$ Axioms]\mbox{}
\begin{enumerate}
\setcounter{enumi}{3}
\item\label{ACF} Any model of $T^\mu$ is an algebraically closed field of characteristic $q$.
\item\label{sol} Given a code $\alpha\in\C_g$ and $\vec{b}$ realizing
$\theta_\alpha(\vec{y})$, one of the following holds:
\begin{enumerate}
\item[a)] $\alpha$ has a pseudo-morley sequence of length
$\mu(\alpha)$ over $\vec{b}$.
\item[b)] Given a realization $\vec{a}$
of $\varphi_\alpha(\vec{x},\vec{b})$ generic over the
model that we are considering, if we paint $\vec{a}$ in black, there
is a code $\beta\in\C_g$ and a pseudo-morley sequence for $\beta$ of
length $\mu(\beta)+1$ in the $\Lans$-structure consisting of the model
and $\vec{a}$ such that there are more than $\mu^*(\beta)$ many
elements of the sequence contained $\{a_1,\dotsc,a_{n_\alpha}\}$.
\end{enumerate}
\end{enumerate}
\end{description}

\begin{note}\label{N:1storder} We discuss here why the above axioms 
are first-order and their meaning. Since our final theory will have
finite Morley rank, it follows from \cite{Mac71} that Axiom (\ref{ACF})
needs to be included.  Axiom (\ref{nosol}) will yield that the types of
$\delta=0$ will become algebraic, and hence of Morley rank 0.

Why is Axiom (\ref{sol}) axiomatizable? In order to encode $\beta$, we
need to determine \emph{a priori} how many variables we will
use. Equivalently, how many $\beta$'s need to be considered. We cannot
use more than $n_\alpha$ variables. On the other hand, we have
$n_\beta$ many variables to consider for each element of the
pseudo-morley sequence, and there are at least $\mu^*(\beta)+1$ many
such members. That is,
\[ (\mu^*(\beta)+1)n_\beta\leq n_\alpha\]

By the finite-to-one condition on $\mu^*$, there are only finitely
many $\beta$'s that satisfy the above inequality, and we are done.

Moreover, it follows from \ref{R:dimcode} that in order to get a
complete theory, we do not need to determine how many realizations of
a code there must be in a model, since we implicitly do so.
\end{note}

\begin{theorem}\label{T:sat=rich}
An $\Lans$-structure is rich if and only if it is an $\omega$-saturated model of $T^\mu$.
\end{theorem}
\begin{proof}
Let $M\models T^\mu$ be $\omega$-saturated. Let $B\leq M$ and $B\leq
A$ be finite sets. We need to find a self-sufficient $B$-copy of $A$
in $M$. Splitting $B\leq A$ into minimal extensions, we are reduced to
the minimal case. We can distinguish four different cases:\\
\noindent
If $B\leq A$ is algebraic, we are done (by Axiom (\ref{ACF})).\\
\noindent
If $B\leq A=B\cup\{\vec{a}\}$ is of type (\ref{black}) (see Lemma
\ref{L:minext}) with $\delta(A/B)=0$, consider the free amalgam
$M'$ of $M$ and $A$ over $B$. Since $M$ is algebraically closed, $M'$
is a good extension of $M$. By (the proof of) Lemma \ref{L:gute-code}
and \ref{L:aux} there is a code $\alpha\in\Cl_g$ and $\vec{b}\in M$
such that $\vec{a}$ is an $M$-generic realization of
$\varphi_\alpha(\vec{x},\vec{b})$. By Axiom (\ref{sol}) $M'$ does not
belong to $\Cl^\mu$. Theorem \ref{T:Amal} implies that $A$ has a
strong embedding over $B$ into $M$.\\
\noindent 
For $1\leq \delta(A/B)\leq p-1$, we need to \emph{approximate} the
extension by extensions of $\delta<\delta(A/B)$ and apply
induction. We know by \ref{L:minext} that $A$ contains no new white
points. Choose some element $a\in A\setminus B$. Since $a$ is
transcendental over $B$, $a^n$ is not in $A$ for large $n$. We can
paint $a^n$ in black and consider $A_n=A\cup\{a^n\}$. It is easy to
check that $B\leq A_n$ is minimal and $\delta(A_n/B)<\delta(A/B)$. The
sequence $A_n/B$ converges (in the space of $\Lan$-types) to the
extension $A_\infty/B$, where $A_\infty=A\cup\{c\}$, with $c$
transcendental over $A$ and black. Clearly $A\leq A_\infty\in\Cl^\mu$,
by \ref{L:member}.  Since there is only a finite number of codes
$\alpha\in\C_g$ for which there could be a pseudo-morley sequence of length longer than $\mu(\alpha)$ in any $A_n$ (bounded only in terms of $|A|$), we have that $A_n$ is
in $\Cl^\mu$ for large $n$. Hence, by induction, we can find
self-sufficient $B$-copies of $A_n$ in $M$ for large $n$. By
saturation of $M$, $A_\infty$ is also self-sufficiently embedable over
$B$. Since $A\leq A_\infty$, we conclude that there is a
self-sufficient $B$-copy of $A$ in $M$.\\
\noindent
For the last case, let $A=B\cup\{a\}$ with $a$ white transcendental
over $B$. Consider for each $n$ the extension
\[B\leq B\cup\{c\}\leq B\cup\{c,c^n\},\]
where $c$ is black transcendental over $B$ and $c^n$ is
white. $B\cup\{c,c^n\}$ belongs to $\Cl^\mu$ by \ref{L:member}. By
the above, that we can realize $B\leq B\cup\{c,c^n\}$
self-sufficiently in $M$. Since these extensions \emph{converge} to
$B\leq B\cup\{c,a\}=A'$ where $c$ and $a$ are algebraically independent
over $B$, we can realize $B\leq A'$ self-sufficiently in
$M$. Since $A\leq A'$, we are done.\\

Suppose now that $M$ is a rich field. We first show that $M$ is
algebraically closed. Let $a\in \acl(M)$. Choose a finite set $B$ in
$M$ such that $a$ is in $\acl(B)$. Taking the closure of $B$ in $M$,
we can assume that $B\leq M$. Paint $a$ in white. It is clear that
$B\cup\{a\}$ is in $\Cl^\mu$ (since $B$ is) and $B\leq B\cup\{a\}$.
By richness, we find a copy of $a$ in $M$ over $B$. This yields
(\ref{ACF}).\\
\noindent
For Axiom \ref{sol}, let $\alpha$ and $\vec{b}$ be as in the statement
such that neither a) nor b) hold. Choose some generic
\emph{black} realization of $\varphi_\alpha(\vec{x},\vec{b})$ over
$M$.  By \ref{C:member}, we have that
$M\cup\vec{a}$ is in $\Cl^\mu$. Choose some finite
set $B\leq M$ containing $\vec{b}$. Again, $B\leq B\cup\vec{a}$, and
by richness, we get a $B$-copy of $\vec{a}$ in $M$, say
$\vec{a}'$. Take now some finite $C\leq M$ containing
$B\cup\vec{a}'$. We have that $C\leq C\cup\vec{a}$. We can iterate and
obtain a pseudo-morley sequence in $M$ for $\alpha$ of arbitrarily
large length. This contradicts that $M$ is in $\Cl^\mu$.\\
\noindent
Now, $M$ is elementarily equivalent to an $\omega$--saturated structure
$M'$, which is by the above a model of $T^\mu$ and therefore rich.  So
$M$ is $\infty$-equivalent to $M'$ and therefore $\omega$--saturated
itself.
\end{proof}

\begin{corollary}\label{C:exclosed}
Let $M$ be an $\Lans$-structure in $\Cl^\mu$. Then $M\models T^\mu$
iff every existential $\Lans(M)$-formula $\phi$ true in some $M\leq N$
with $N\models\textrm{``Universal Axioms of $T^\mu{}$''}$ holds also
in $M$.
\end{corollary}
\begin{proof}
Let $M$ be a model of $T^\mu$, $N$ and $\phi$ as above. We can assume
$M$ and $N$ are saturated. Let $B\leq M$ contain all parameters in
$\phi$. Choose some $B\subset A\leq N$ containing a realization of
$\phi$. Since $B\leq A$ and $M$ is rich, we can embed $A$ in $M$ over
$B$. Hence, we have a solution for $\phi$ in $M$.

If $M$ is existentially closed among self-sufficient extensions, it
satisfies Axiom (\ref{ACF}), since $\acl(M)$ (new elements painted in
white) is a self-sufficient extension. If $M$ does not satisfy Axiom
(\ref{sol}), there is an $\alpha\in\C_g$, $\vec{b}\in M$ and a black
$M$-generic solution $\vec{a}$ of $\varphi_\alpha(\vec{x},\vec{b})$
such that $M\leq M\cup\vec{a}$ is in $\Cl^\mu$. Considering finite
sets $C\leq M$ containing $\vec{b}$ and using the existential
closedness of $M$, we can find infinitely many disjoint black
realizations of $\varphi_\alpha(\vec{x},\vec{b})$. (We may assume that
$\varphi_\alpha$ is quantifier free.) This contradicts
that $M\in \Cl^\mu$.
\end{proof}

\begin{corollary}\label{C:vollstaendig}
  The theory $T^\mu$ is complete. 

  Two tuples $\vec{a}$ and $\vec{b}$ in two
  models $M$ and $N$ of $T^\mu$ have the same $\Lans$-type iff there is an
  isomorphism $f:\cl(\vec{a})\to\cl(\vec{b})$ which maps $\vec{a}$ to
  $\vec{b}$. 

  An extension $M\subset N$ of models of $T^\mu$ is elementary iff $M$
  is self-sufficient in $N$\footnote{By an observation of M. Hils
  $T^\mu$ is model complete, i.e.\ \emph{all} extensions of models of
  $T^\mu$ are self-sufficient. See Remark \ref{hils} for a proof.}
\end{corollary}
\begin{proof}
$T^\mu$ is complete, since any two countable saturated
models are elementarily equivalent, by richness.

Consider two models $M$, $N$ of $T^\mu$. If $M$ is subset of $N$,
but no self-sufficient, there is a finite $A\leq M$ and a tuple
$\vec{b}\in N$ sucht that $\delta(\vec{b}/A)<0$. Responsible is a
finite part of the $\Lans$-type of $\vec{b}$ over $A$. So $M\prec N$
would imply the existence of an $\vec{a}\in M$ with
$\delta(\vec{a}/A)<0$, which is not possible.

If $\vec{a}$ and $\vec{b}$ have the same $\Lans$--type, it is easy to
see that the map $\vec{a}\mapsto\vec{b}$ extends to an isomorphism
$f:\cl(\vec{a})\to\cl(\vec{b})$. Conversely, let  $f$ be given. Choose
rich extensions $M\prec M'$ and $N\prec N'$. We know that $M$ and $N$
are self-sufficient in these extension and therefore also
$cl(\vec{a})$ and $\cl(\vec{b})$. Since isomorphisms between finite
self--sufficent subsets of $M'$ and $N'$ have the back-and-forth
property, $f$ is elementary map.

Finally assume that $M\leq N$. Since $\cl_M$ is the restriction of
$\cl_N$ to $M$, all finite tuples $\vec{a}\in M$ have the same
$\Lans$-type in $M$ as in $N$, i.e. $M\prec N$.
\end{proof}

\section{Computing ranks}\label{S:ranks}

In this section, we compute the Morley rank of types in
$T^\mu$. In order to avoid confusion, we will  denote it by $\Mrs$, since we work with $\Lans$-types $\tps(\vec{a}/B)$. We work inside a sufficiently saturated
model $M$ of $T^\mu$.

\begin{lemma}\label{L:fteMR}
$T^\mu$ has finite Morley rank.
\end{lemma}
\begin{proof}
  It is clear that $\cl(A)$ is contained in $\acls(A)$. This implies
  $$\Mrs(\vec{a}/C)=\Mrs(\cl(C\vec{a})/\cl(C)).$$ 
  
  So it is enough to
  compute Morley ranks $\Mrs(A/B)$, where $B\leq A\leq M$ and
  $A\setminus B$ is finite. We will show that the rank is bounded by a
  function of $\delta(A/B)$.

  We prove first that $\delta(A/B)=0$ implies that $A$ is algebraic in
  $M$ over $B$, i.e.\ \linebreak $\Mrs(A/B)=0$. For this we may
  assume that $A/B$ is minimal. If $A$ has a new white element, then
  $A/B$ is algebraic (in the field sense). Otherwise, $A\setminus B$
  contains only black points and we may assume -- after adding algebraic elements to $B$
  -- that $A/B$ is good. By \ref{L:gute-code}, $A\setminus B$ is
  enumerated by a generic solution $\vec{a}$ of a code $\alpha\in\C_g$
  over some $b\in\dcl(B)$. By \ref{L:aux} any sequence $A_i$ of
  different conjugates of $A$ over $B$ in $M$ yields a sequence of
  $B$--generic realizations of $\varphi_\alpha(\vec{x},\vec{b})$ which
  is (in $\K$) independent over $B$. So the sequence is a
  pseudo-morley sequence of $\alpha$ over $\vec{b}$ and cannot be
  longer that $\mu(\alpha)$. This proves that $A/B$ is algebraic in
  $M$.\footnote{Note that the number of conjugates of $\vec{a}$ over
  $B$ is bounded by $n_\alpha!\cdot\mu(\alpha)$.}\\

  \noindent Now, assume that $\delta(A/B)=d>0$ and that $\Mrs(A'/B')\leq
  r$ for all $B'\leq A'\leq M$ with $\delta(A'/B')<d$. The above
  case shows that we may also assume that for all $C\leq A$ with
  $B\subsetneq C\subsetneq A$ we have $\delta(C/B)>0$ and
  $\delta(A/C)>0$. We distinguish three cases:\\
  
  \noindent Case 1: $A/B$ is not minimal. Then we find $B\leq C\leq A$
  with $\delta(C/B)< d$ and $\delta(A/C)< d$.
  Enumerate $C\setminus B$ by $\vec{c}$ and $A\setminus C$ by
  $\vec{a}$ and choose $\Lans$-formulas $\rho(\vec{z})$ and
  $\chi(\vec{z},\vec{x})$ over $B$ which are satisfied by $\vec{c}$
  and $\vec{c}\vec{a}$ and which imply $\delta(\vec{z}/B)<d$ and
  $\delta(\vec{x}/B\vec{z})<d$. By above inductive assumption, we have
  $\Mrs(\vec{c}'/B)\leq r$ and
  $\Mrs(\vec{a}'/B\vec{c}')\leq r$ for all realization $\vec{c}'\vec{a}'$
  of $\rho(\vec{z})\land\chi(\vec{z},\vec{x})$.
  Now we can apply Erimbetov's inequalities \cite{mE75} and obtain
  $\Mrs(\rho(\vec{z})\land\chi(\vec{z},\vec{x}))\leq r(r+1)$. Hence
  $\Mrs(A/B)\leq r(r+1)$.\\

  \noindent Case 2: $A/B$ is minimal and $1\leq d\leq p-1$. We fix an
  enumeration $\vec{a}$ of $A\setminus B$.  We may again assume that
  $p=\tp(\vec{a}/B)$ is stationary. Choose an $\Lan$-formula
  $\phi(\vec{x})$ in $p$ of the same Morley rank and of degree $1$
  which satifies \ref{D:code} (\ref{D:code_algebraisch}).  It follows from Lemma \ref{L:aux} that for every black realization
  $\vec{a}'$ of $\varphi$, only two possibilities may occur: \\ 
  Either we have
  $\delta(\vec{a}'/B)<d$, which implies $\delta(\cl(B\vec{a}')/B)<d$
  and $\Mrs(\vec{a}'/B)=\Mrs(\cl(\vec{a}')/B) \leq r$, or we have
  $\delta(\vec{a}'/B)=d$, which implies $\tp(\vec{a}'/B)=p$. If
  $B\vec{a}'$ is not self-sufficient in $M$, we conlude again
  $\Mrs(\cl(\vec{a}')/B) \leq r$. Otherwise, by \ref{C:vollstaendig}, we have that $\tps(\vec{a}'/B)=\tps(\vec{a}/B)$. This shows
  that the $\Lans$--type of $\vec{a}$ over $B$ is not an accumulation
  point of types over $B$ of rank bigger than $r$.  If $B$ were an
  $\omega$-saturated elementary substructure of $M$, we could
  conclude that $\Mrs(A/B)\leq r+1$. Hence, consider any $\omega$-saturated
  elementary substructure $N$ which contains $B$.  Since $A$ is either
  contained in $N$ or intersects $N$ in $B$, we have
  $\delta(A/N)\leq\delta(A/B)$, which implies $\Mrs(A/N)\leq
  r+1$. Since $N$ was arbitrary, it follows $\Mrs(A/B)\leq r+1$.\\

  \noindent Case 3: $A/B$ is minimal and $d=p$.  Then $A=B\cup\{a\}$
  with $a$ white and transcendental. All $1$-Types different from
  $\tps(a/B)$ have rank $\leq r$. So, if $B$ were an
  $\omega$--saturated elementary model, we could conclude
  $\Mrs(a/B)\leq r+1$. By the same argument as above, we show that the claim holds.
\end{proof}

For any set of parameters $C$ and any finite tuple $\vec{a}$ we define
$$\dd(\vec{a}/C)=\delta(\cl(C\vec{a})/\cl(C)).$$

\begin{theorem}\label{RM=d}
$\Mrs(\vec{a}/C)=\dd(\vec{a}/C)$.
\end{theorem}

\begin{proof}
Since $M$ is a field of finite Morley rank, Morley rank satisfies
the Lascar inequalities by a result of Lascar \cite{dL85}. It follows
now from the proof of \ref{L:fteMR} and the additive character of
$\delta$ that $\Mrs\leq d$.

For the other inequality, let $B\leq A\leq M$ be minimal, and
$\delta(A/B)=d>0$. It follows from the proof of \ref{T:sat=rich} (and
from \ref{C:vollstaendig})
there is a sequence of extensions $B\leq A_n\leq M$, such that 
$\delta(A_n/B)=d-1$ and
$$\lim_{n\to\infty}\tps(A_n/B)=\tps(A/B).$$ Since $\Mrs(A_n/B)\geq
d-1$ by induction, we have $\Mrs(A/B)\geq d$.
\end{proof}

\begin{proof}[Proof of the Main Theorem:]
  For any $a$ in $M$, we have:
  $$\dd(a)\leq\delta(a)=
  \begin{cases}
    p & \text{if $a$ is white}\\
    p-1 & 
    \text{if $a$ is black}
  \end{cases}$$
  This shows $\Mrs(T^\mu)\leq p$ and $\Mrs(N)\leq p-1$.
  On the other hand the structures $\{a\}$, $\{b\}$  with $a$ white,
  transcendental and with $b$ black
  transcendental are both in $\Cl^\mu$. So we find them as
  self-sufficient subsets of $M$. Then $\dd(a)=\delta(a)=p$
  and $\dd(b)=\delta(b)=p-1$. This proves the result.
\end{proof}

\begin{example}\label{E:sumblacks}
It was observed in Theorem 18 \cite{BH01} that every generic white point is the sum of
$p$ independent black points of Morley rank $1$. We want to give a simpler proof of this fact. 

Let $a_1,\dotsc,a_p$ be generic independent elements of $\K$.
It follows trivially from \ref{L:member} that
the black $A_i=\{a_i^r\}_{1\leq r\leq p-1}$
and $\bigcup\limits_{i=1}^{p} A_i$ with all elements painted in black
belong to $\Cl^\mu$.

Take now $a=a_1+\dotsb+a_p$ painted in white. Again by \ref{L:member},
$\bigcup\limits_{i=1}^{p} A_i \cup\{a\}$ belongs to $\Cl^\mu$. We
may assume that $\bigcup\limits_{i=1}^{p} A_i \cup\{a\}\leq M$. 

Since $\{a\}\leq M$, we have $d(a)=p$. Hence, $a$ is an generic white element.
Since $\cl(a_i)=A_i\leq M$, we have $d(a_i)=\delta(A_i)=1$. So each $a_i$ has Morley rank $1$.
\end{example}
\begin{remark}[added January 23, 2005]\label{hils}
  Martin Hils made the following oberservation: A theory of fields of
  finite Morley rank is $\aleph_1$--categorical. Since $T^\mu$ is
  $\forall\exists$--axiomatizable, a theorem of Lindstr\"om implies
  that \emph{ $T^\mu$ is model complete.}  
  
  Let us give a direct proof. Assume that $M$ is a model of $T^\mu$
  and $N$ an extension which belongs to $\Cl^\mu$. We want to show
  that $M$ is self--sufficient in $N$. We may assume that
  $\delta(N/M)<0$ and that $N$ is minimal with this property. Then
  $N=N'\cup\{a\}$, where $M\leq N'$ with $\delta(N'/M)=0$, and $a$ is
  black and algebraic over $N'$ in the field sense. Choose a rich
  field $N'\leq N''$. Then $N''$ is an elementary extension of $M$. On
  the other hand, we have $\Mrs(N'/M)=\dd(N'/M)=\delta(N'/M)=0$, so
  $N'$ is, in $N''$, algebraic over $M$ and therefore contained in
  $M$. This implies that $a$ is field-algebraic over $M$, which is
  impossible, since $M$ is an algebraically closed field.
\end{remark}


\begin{thebibliography}{99}

\bibitem{BH00} J. Baldwin, K. Holland, \emph{Constructing $\omega$-stable structures: rank 2 fields}, J. Symb. Logic, {\bf 65}, $\mathrm{n}^\mathrm{o}$ 1, 371--391, (2000).

\bibitem{BH01} J. Baldwin, K. Holland, \emph{Constructing $\omega$-stable structures: Rank k fields}, Notre Dame J. of Formal Logic, {\bf 44}, 139--147, (2004).


\bibitem{BZ04} A. Baudisch, M. Ziegler, \emph{Hrushovskis Fusion}, preprint (2004).

\bibitem{Bu04} Juan R. Bueno, \emph{Le carr\'e de l'\'egalit\'e collaps\'e}, preprint (2004).

\bibitem{mE75} M. Erimbetov, \emph{Complete theories with 1-cardinal formulas}, Algebra i Logika, {\bf 14}, $\mathrm{n}^\mathrm{o}$ 3, 245--257 (1975).

\bibitem{Hr92} E. Hrushovski, \emph{Strongly minimal expansions of algebraically closed fields}, Israel J. Math, {\bf 79}, 129--151, (1992).

\bibitem{Hr93} E. Hrushovski, \emph{A new strongly minimal set}, Annals of Pure and Applied Logic, {\bf 62}, 147--166, (1993).

\bibitem{HZ93} E. Hrushovski, B. Zil'ber, \emph{Zariski Geometries}, Bull. Amer. Math. Soc. {\bf 28}, 315--323, (1993).

\bibitem{dL85} D. Lascar, \emph{Les groupes $\omega$-stables de rang fini}, Trans. Amer. Math. Soc., {\bf 292}, 451--462, (1985).

\bibitem{Mac71} A. Macintyre, \emph{On $\omega_{1}$-categorical fields}, Fund. Math.,{\bf 71}, 1-25 (1971).

\bibitem{Po99} B. Poizat, \emph{Le carr\'e de l'\'egalit\'e }, J. Symb. Logic, {\bf 64}, $\mathrm{n}^\mathrm{o}$ 3, 1338--1355, (1999).



\bibitem{BMZ04} A. Baudisch, A. Martin-Pizarro, M. Ziegler,
\emph{Hrushovskis Fusion}, preprint (2004).

\end{thebibliography}
\end{document}